\documentclass[english,12pt]{article}
\usepackage[T1]{fontenc}
\usepackage{babel}
\usepackage{amsmath,amsthm}
\usepackage{amssymb}
\usepackage[all]{xy}

\newcommand{\reel}{\mathbb{R}}
\newcommand{\comp}{\mathbb{C}}
\newcommand{\rat}{\mathbb{Q}}
\newcommand{\rel}{\mathbb{Z}}

\newcommand{\ra}{\rightarrow}

\newtheorem{thm}{Theorem}

\newtheorem{prop}[thm]{Proposition}
\newtheorem{lem}[thm]{Lemma}

\title{Conjugate varieties with distinct real cohomology algebras}
\author{Fran\c{c}ois Charles\\ \'Ecole Normale Sup\'erieure, Paris}

\begin{document}

\maketitle

\begin{abstract}
Using constructions of Voisin, we exhibit a smooth projective variety defined over a number field $K$ and two complex embeddings of $K$, such that the two complex manifolds induced by these embeddings have non isomorphic cohomology algebras with real coefficients. This contrasts with the fact that the cohomology algebras with l-adic coefficients are canonically isomorphic for any prime number l, and answers a question of Grothendieck.
\end{abstract}


\section{Introduction}

Let $X$ be an algebraic variety defined over an algebraically closed field $K$ of characteristic $0$. If $l$ is a prime number, the $l$-adic cohomology of $X$, $H^*(X, \rat_l)$, is a graded algebra over $\rat_l$. Its definition as an inverse limit of \'etale cohomology groups shows that it does not depend on the structural map $X\ra\mathrm{Spec}\,K$.

Now suppose $K$ is a finitely generated extension of $\rat$. We can consider the $l$-adic cohomology of $X_{\bar{K}}$, where $\bar{K}$ is an algebraic closure of $K$. This does not depend on the choice of the algebraic closure. Moreover, if $\bar{L}$ is any algebraically closed field containing $K$, the proper base change theorem shows that the $l$-adic cohomology of $X_{\bar{L}}$ is canonically isomorphic to that of $X_{\bar{K}}$. In particular, $H^*(X_{\comp}, \rat_l)$ is canonically defined and does not depend on an embedding $K\hookrightarrow \comp$. Let us choose such an embedding, and assume from now on that $X$ is smooth over $K$. Artin's comparison theorem of \cite{SGA4}, Exp. XI, shows that $H^*(X_{\comp}, \rat_l)$ is canonically isomorphic to the Betti cohomology of $X_{\comp}^{an}$, the underlying complex manifold of $X_{\comp}$, with coefficients in $\rat_l$. As a consequence, the latter does not depend on the embedding $K\hookrightarrow \comp$.

The preceding discussion shows that the cohomology algebra with coefficients in $\rat_l$ of the complex manifold $X_{\comp}^{an}$ does not vary under automorphisms of $\comp$. In other words, it does only depend on the abstract scheme $X_{\comp}$ and not on its map to $\mathrm{Spec}\,\comp$.

The topology of a complex variety can nonetheless vary under automorphisms of $\comp$. Indeed, Serre constructs in \cite{S} two conjugate complex smooth projective varieties with different fundamental groups. Note however that the profinite completions of those are canonically isomorphic, due to Grothendieck's theory of the algebraic fundamental group. In particular, they do not have the same homotopy type. 

Other constructions of conjugate varieties which are not homeomorphic can be found in \cite{A}, and more recently in \cite{C} and \cite{Shi} (the last article actually considers open varieties). See also \cite{ACC} for related constructions. Nevertheless, the arguments leading to the constructions of the previous examples all make use of the integral homotopy type, by actually considering fundamental groups or Betti cohomology with coefficients in $\rel$.

This leads naturally to the following question, asked in \cite{R} : 
\begin{center}
Do there exist conjugate varieties with different rational homotopy type ?
\end{center}
Very similarly, it had already been asked by Grothendieck (Montr\'eal, july 1970) whether there exist conjugate varieties with distinct cohomology algebras with rational coefficients.

In this paper, we answer positively these questions. We actually show the following, which is stronger : 
\begin{thm}
There exist smooth projective conjugate varieties whose real cohomology algebras are not isomorphic.
\end{thm}

To construct the example, we use the methods of \cite{V}, where Voisin shows how to use the cohomology algebra of some varieties to recover their endomorphism rings. As in \cite{S}, our example is built out of abelian varieties with complex multiplication.

\paragraph{Acknowledgements.}
I would like to thank Claire Voisin for sharing her ideas with me, and for her many very helpful remarks. This paper owes her a lot. I would also like to thank Johan de Jong for many useful discussions. This work was completed during a stay at Columbia University.

\section{Statement of the theorem}

Let $k$ and $k'$ be two different imaginary quadratic subfields of $\comp$ \footnote{For simplicity, but at the loss of some functoriality, we fix complex embeddings of the quadratic fields.}, and let $E$ (resp. $E'$) be a complex elliptic curve with complex multiplication by $\mathcal{O}_k$ (resp. $\mathcal{O}_{k'}$). Let $A$ be the product of $E$ and $E'$. 

Suppose we are given polarizations $\phi$ and $\phi'$ of $E$ and $E'$, that is, numerical equivalence classes of very ample line bundles. Those give a polarization $\psi$ of $A$ with $\psi=\phi\oplus\phi'$, which comes from some projective embedding $i:A\hookrightarrow \mathbb{P}^N$. The idea of the example is to construct a variety whose cohomology ring encodes the endomorphism ring of $A$ and contains a distinguished line related to the polarization. This will be achieved by blowing-up some special subvarieties of $A\times A\times \mathbb{P}^N$.

Let $a$ and $a'$ be  imaginary elements of $\mathcal{O}_k$ and $\mathcal{O}_{k'}$ which generate the fields $k$ and $k'$ respectively.  We will denote by $f$ (resp. $f'$) the endomorphism $a\times 0$ (resp. $0\times a'$) of $A=E\times E'$ -- we will sometimes still denote by $f$ (resp. $f'$) the corresponding endomorphism of $E$ (resp. $E'$). The induced homomorphisms on the first cohomology group of $A$ have eigenvalues $0$ and $a, -a$ (resp. $a', -a'$).

Let $x, y, z$ and $t$ be points of $\mathbb{P}^N$ and $u$ be a point of $A$. Let us consider the following smooth subvarieties of $A\times A\times \mathbb{P}^N$ : 
\begin{displaymath}
Z_1=A\times 0\times x,\, Z_2=\Gamma_{\mathrm{Id}_A}\times y, \,Z_3=\Gamma_{f}\times z, \,Z_4=\Gamma_{f'}\times t,\, Z_5=u\times \Gamma_i,
\end{displaymath}
where $\Gamma$ stands for the graph of a morphism. For a generic choice of $x, y, z, t, u$, those subvarieties are pairwise disjoint. Let $X$ be the blow-up of $A\times A\times\mathbb{P}^N$ along those subvarieties. This is a smooth complex projective variety\footnote{Since abelian varieties with complex multiplication are defined over number fields, we can even choose the polarizations and the points adequately so that $X$ is also defined over a number field.}.

\bigskip

For any smooth complex variety $V$, one can consider the underlying complex manifold $V^{an}$ and its real cohomology algebra $H^*(V^{an}, \reel)$. By an abuse of notation, we will denote the latter algebra by $H^*(V, \reel)$, remembering that it depends on the usual topology on $V$. For a scheme $Y$ over a field $K$, and $\sigma$ an automorphism of $K$, let $Y^{\sigma}$ denote the scheme $Y\otimes_{K, \sigma} K$. If $Y$ is smooth (resp. polarized), then so is $Y^{\sigma}$. Our theorem is the following.

\begin{thm}\label{th}
Let $X$ be constructed as above, and let $\sigma$ be an automorphism of $\comp$ which acts trivially on one of the fields $k$ and $k'$, but not on the other. Then the real cohomology algebras $H^*(X, \reel)$ and $H^*(X^{\sigma}, \reel)$ are not isomorphic.
\end{thm}

\section{Proof of the theorem}

Theorem \ref{th} will be obtained as a consequence of propositions \ref{CM} and \ref{main}, which will be stated and proved in the next subsections.

\subsection{Some linear algebra}\label{la}

Let $\sigma$ be an automorphism of the field $\comp$. In this section, we describe some linear objects attached to the polarized abelian variety $(A^{\sigma}, \psi^{\sigma})$ and study how they vary under the action of automorphisms of $\comp$. 
What we would like to do is to recover the CM-type of $E^{\sigma}$ and $E'^{\sigma}$ from part of the cohomology of $X^{\sigma}$. Actually, we will only be able to compare, in some sense, those CM-types, which will be enough for our purpose. This is the reason why we have to work with two elliptic curves.

\bigskip

Let $\sigma$ be an automorphism of $\comp$. There is a canonical isomorphism of abstract schemes from $A^{\sigma}$ to $A$. Though it is by no means defined over $\comp$, it still induces an isomorphism between the endomorphism rings of the complex varieties $A^{\sigma}$ and $A$. As a consequence, there is a canonical action of $k\times k'$ on $A^{\sigma}$. The real vector space $H^1(A^{\sigma}, \reel)$ thus becomes a free rank $1$ $\comp\times \comp$-module, as $k\otimes_{\rat}\reel$ and $k'\otimes_{\rat}\reel$ are canonically isomorphic to $\comp$ (recall we chose complex embeddings of $k$ and $k'$). From the embedding $i^{\sigma} : A^{\sigma}\hookrightarrow \mathbb{P}^N$, we get a homomorphism $i^{\sigma *} : H^2(\mathbb{P}^N, \reel)\ra \bigwedge^2_{\reel} H^1(A^{\sigma}, \reel)$.

As a consequence, with each $\sigma$ comes a free rank $1$ $\comp\times \comp$-module 
$$V=H^1(A^{\sigma}, \reel),$$
a $1$-dimensional $\reel$-vector space 
$$L=H^2(\mathbb{P}^N, \reel)$$ 
and a nonzero homomorphism 
$$\mu=i^{\sigma *} : L\ra \bigwedge^2_{\reel} V.$$

Let $\mathcal{L}$ be the set of isomorphism classes of such triples $(V, L, \mu)$, with the obvious notion of morphism. The preceding description gives us a map

$$l : \mathrm{Aut}(\comp) \ra \mathcal{L}.$$

\begin{prop}\label{CM}
Let $\sigma\in\mathrm{Aut}(\comp)$ act trivially on one of the fields $k$ and $k'$, but not on the other. Then $l(\sigma)\neq l(\mathrm{Id}_{\comp})$.
\end{prop}

\begin{proof}
Let $V$ be a free rank $1$ $\comp\times \comp$-module. Using the idempotents of $\comp\times\comp$, we get a canonical splitting of $V$ as a direct sum of two complex vector spaces $V_1$ and $V_2$ of rank $1$, such that $0\times\comp$ acts trivially on $V_1$ and $\comp\times 0$ acts trivially on $V_2$. The action of $\comp\times 0$ on $V_1$ and of $0\times\comp$ on $V_2$ endows those real vector spaces with a complex structure, so the underlying real vector spaces of $V_1$ and $V_2$ are canonically oriented. Let us call those complex structures and the induced orientations the \emph{algebraic} ones.

Let $L$ be a real vector space of rank $1$. Any homomorphism $\mu$ of real vector spaces from $L$ to $\bigwedge^2_{\reel} V$ canonically induces homomorphisms from $L$ to $\bigwedge^2_{\reel} V_1$ and $\bigwedge^2_{\reel} V_2$. Suppose that those are isomorphisms -- this is the case for the triples in the image of $l$. We get an isomorphism between $\bigwedge^2_{\reel} V_1$ and $\bigwedge^2_{\reel} V_2$, which may or may not respect the algebraic orientation. Let us define the \emph{sign} of the triple $(V, L, \mu)$ to be $1$ or $-1$ according to wether it is the case or not. The sign only depends on the isomorphism class of the triple.

\smallskip

In the setting of the proposition, the sign of a triple is easy to compute. Let us indeed choose an automorphism $\sigma$ of $\comp$. The aforementioned splitting of $H^1(A^{\sigma}, \reel)$ corresponds to the splitting 
$$H^1(A^{\sigma}, \reel)=H^1(E^{\sigma}, \reel)\oplus H^1(E'^{\sigma}, \reel).$$

Aside from the complex structure induced by complex multiplication, the space $H^1(E^{\sigma}, \reel)$ has a complex structure induced by its identification with the cotangent space at $0$ to the complex manifold $E^{\sigma}$. It does not have to agree with the one previously defined through complex multiplication. The same construction works with $E'$. Let us call those complex structures and the induced orientations \emph{transcendental}. Now let $h\in H^2(\mathbb{P}^N, \reel)$ be the class of a hyperplane. The homomorphisms $H^2(\mathbb{P}^N, \reel)\ra \bigwedge^2_{\reel} H^1(E^{\sigma}, \reel)$ and $H^2(\mathbb{P}^N, \reel)\ra \bigwedge^2_{\reel} H^1(E'^{\sigma}, \reel)$ defined earlier both send $h$ to elements which are positive with respect to the transcendental orientation. 

As a consequence, the isomorphism $\bigwedge^2_{\reel} H^1(E^{\sigma}, \reel)\ra \bigwedge^2_{\reel} H^1(E'^{\sigma}, \reel)$ induced by the polarization respects the transcendental orientations. From this remark, it results that the sign of the triple $l(\sigma)$ is $1$ if and only if the algebraic and transcendental orientations either coincide on both $H^1(E^{\sigma}, \reel)$ and $H^1(E'^{\sigma}, \reel)$, or if they differ on both those spaces.

\smallskip

Recall that a complex structure on a real vector space $V$ is given by an automorphism $I$ of $V$ such that $I^2=-\mathrm{Id}_V$. Giving $I$ is in turn equivalent to giving a splitting of the complex vector space 
$$V\otimes\comp=V^{1,0}\oplus V^{0,1}$$
such that $V^{1,0}$ and $V^{0,1}$ are complex conjugate of each
other. Those spaces are the eigenspaces of $I$ for the eigenvalues $i$ and
$-i$ respectively.

Now let $I^{\sigma}_{alg}$ and $I^{\sigma}_{tr}$ be the
automorphisms of $H^1(E^{\sigma}, \reel)$ corresponding to the algebraic
and transcendental complex structures. Since the action of complex
multiplication on $H^1(E^{\sigma}, \reel)$ is $\comp$-linear with respect
to the transcendental complex structure (indeed, morphisms of smooth
complex algebraic
varieties are holomorphic), $I^{\sigma}_{alg}$ and $I^{\sigma}_{tr}$
commute, which implies, as their eigenspaces are one-dimensional and they
have $i$ and $-i$ has eigenvalues, that they are either equal or opposite
to each other.

The splitting of $H^1(E^{\sigma}, \reel)$ corresponding to the 
transcendental complex structure is well-known, as it corresponds to the 
Hodge decomposition $$H^1(E^{\sigma}, 
\comp)=H^0(E^{\sigma},\Omega_{E^{\sigma}})\oplus 
H^1(E^{\sigma},\mathcal{O}_{E^{\sigma}}),$$ with $I^{\sigma}_{tr}$ 
acting as $i$ on the first summand and as $-i$ on the second. Therefore, 
we just have to investigate the action of complex multiplication on 
$H^0(E^{\sigma},\Omega_{E^{\sigma}})$.

We have an obvious isomorphism of one-dimensional complex vector spaces 
$$H^0(E,\Omega_{E})\ra H^0(E^{\sigma},\Omega_{E^{\sigma}}), 
\omega\mapsto \omega^{\sigma}$$ given by pullback of differential forms 
by the isomorphism of abstract schemes $E^{\sigma}\ra E$. This 
isomorphism is $\sigma$-linear, that is, it sends $\lambda\omega$ to 
$\sigma(\lambda)\omega^{\sigma}$. Let $f$ be the endomorphism of $E$ we 
chose 
earlier. It generates its complex multiplication, and corresponds to a 
certain imaginary $a\in \mathcal{O}_k\subset \comp$. For any global 
holomorphic one-form $\omega$ on $E$, we have 
$$f^{\sigma*}\omega^{\sigma}=(f^*\omega)^{\sigma}.$$ The morphisms $f^*$ 
and $f^{\sigma*}$ act as scalars on $H^0(E,\Omega_{E})$ and 
$H^0(E^{\sigma},\Omega_{E^{\sigma}})$ respectively. Let $f^*$ act by 
multiplication by $\lambda$. The complex number $\lambda$ may be either 
$a$ or $\bar{a}=-a$, and it is equal to $a$ if and only if the algebraic 
and transcendental complex structures on $H^1(E, \reel)$ coincide. But 
$$f^{\sigma*}\omega^{\sigma}=(f^*\omega)^{\sigma}=(\lambda\omega)^{\sigma}=\sigma(\lambda)\omega^{\sigma}.$$

This proves that if $\sigma$ is an automorphism of $\comp$ fixing $a$, hence $k$, then the transcendental and algebraic complex structures on $H^1(E^{\sigma}, \reel)$ coincide if and only if they do on $H^1(E, \reel)$. The same goes for the algebraic and transcendental orientations. On the other hand, if $\sigma$ acts as complex conjugation on $k$, the transcendental and algebraic orientations on $H^1(E^{\sigma}, \reel)$ coincide if and only if they don't on $H^1(E, \reel)$.

Since the same goes for $E'^{\sigma}$ and $k'$, the preceding discussion shows that the sign of $l(\sigma)$ is the same as the sign of $l(\mathrm{Id}_{\comp})$ if and only if $\sigma$ acts either trivially on $k$ and $k'$ or by complex conjugation on both. This concludes.
\end{proof}

\paragraph{Remark.}
Using the Hasse principle and the main theorem of complex multiplication, one can prove that the image of $l$ has exactly two elements, and that $l(\sigma)=l(\mathrm{Id}_{\comp})$ if either $\sigma$ acts trivially on both $k$ and $k'$ or if its acts by complex conjugation on both. This is still true if we replace elliptic curves by abelian varieties (assuming the polarizations are compatible with the complex multiplication).

\subsection{Analysis of the cohomology algebra}

The goal of this section is to prove the following.

\begin{prop}\label{main}
The variety $X$ being defined as in the previous section, let $\mathcal{C}$ be the set of isomorphism classes of real graded algebras of the form $H^*(X^{\sigma}, \reel)$, where $\sigma$ is an automorphism of $\comp$, and let $c$ be the map
$$\mathrm{Aut}(\comp)\ra \mathcal{C}, \sigma \mapsto H^*(X^{\sigma}, \reel).$$
Then the map $l:\mathrm{Aut}(\comp)\ra\mathcal{L}$ defined in \ref{la} factors through $c$.
\end{prop}
In other words, given the real cohomology algebra of $X^{\sigma}$, one can recover the linear algebra data described previously. 

\bigskip

Let us fix some notations. Let $\tau : X\ra A\times A\times \mathbb{P}^N$ be the blowing-down map and $\pi: X\ra A\times A$ be the composite of $\tau$ with the projection on the first two factors. While $Z_1, \ldots, Z_5$ are the centers of the blow-up, let $D_1, \ldots, D_5$ denote the corresponding exceptional divisors of $X$ and for $k$ between $1$ and $5$, let $$j_k : Z_k\hookrightarrow A\times A\times \mathbb{P}^N, \,\widetilde{j}_k : D_k\hookrightarrow X, \tau_k=\tau_{|D_k}$$ be the canonical morphisms. For a subvariety $Z$ of $X$, let $[Z]$ denote its cohomology class. Let $h\in H^2(X, \reel)$ be induced by the cohomology class of a hyperplane in $\mathbb{P}^N$ via the natural morphism $X\ra \mathbb{P}^N$.
We will use the same notations, with a superscript $\sigma$, for the corresponding objects of $X^{\sigma}$.

\bigskip

Proposition \ref{main} is a straightforward consequence of the following two statements.

\begin{prop}\label{morph}
Let $\mathcal{C}'$ be the set of isomorphism classes of 7-uples $(H, L_1, \ldots, L_6)$, where $H$ is an element of $\mathcal{C}$ and $L_1, \dots, L_6$ are lines in $H$. Let $c'$ be the map
$$\mathrm{Aut}(\comp)\ra \mathcal{C}', \sigma \mapsto (H^*(X^{\sigma}, \reel), \reel[D_1^{\sigma}], \ldots, \reel[D_5^{\sigma}], \reel h^{\sigma}).$$
If $\sigma$ and $\sigma'$ are automorphisms of $\comp$, then $c(\sigma)=c(\sigma')$ if and only if $c'(\sigma)=c'(\sigma')$.
\end{prop}

\begin{prop}\label{recover}
The map $l$ factors through $c'$.
\end{prop}

Before going through the proofs, let us state a lemma for future reference. This is the analogue of the computations made in the proof of theorem 3 of \cite{V}, and, as in Voisin's paper, is the key to extracting information from the algebra structure on cohomology spaces. We give a proof for the reader's convenience.

\smallskip

Let $p_1$ (resp. $p_2$) be the restriction map from $H^1(A\times A, \reel)$ to $H^1(A, \reel)$ induced by the inclusion of the first (resp. second) factor. Let $q_1$ (resp. $q_2$) be the restriction map from $H^2(A\times\mathbb{P}^N, \reel)$ to $H^2(A, \reel)$ (resp. $H^2(\mathbb{P}^N, \reel)$) induced by the inclusion of the first (resp. second) factor. Using pullback by $\pi$, we can consider $H^*(A\times A, \reel)$ as a subspace of $H^*(X, \reel)$. Note that $H^1(A\times A, \reel)=H^1(X, \reel)$.

Let us fix a nonzero cohomology class $\alpha\in H^2(A\times A, \reel)\subset H^2(X, \reel)$. For $k$ between $1$ and $5$, consider cup-product with $[D_k]$. It gives a homomorphism $H^1(A\times A, \reel)=H^1(X, \reel)\ra H^{3}(X, \reel)$. We get similar homomorphisms  by taking cup-product with $h$ or $\alpha$.

\begin{lem}\label{cup}
The images of those homomorphisms $\cup [D_1], \ldots, \cup [D_5],\cup h,$ and $\cup \alpha$, are in direct sum. 
Furthermore, Ker$(\cup \alpha)$ is at most $2$-dimensional and

\begin{itemize}
\item $\mathrm{Ker} (\cup h)=0, $
\item $\mathrm{Ker} (\cup [D_1])=\mathrm{Ker}(p_1), $
\item $\mathrm{Ker} (\cup [D_2])=\mathrm{Ker}(p_1+p_2), $
\item $\mathrm{Ker} (\cup [D_3])=\mathrm{Ker}(p_1+f^*\,p_2), $
\item $\mathrm{Ker} (\cup [D_4])=\mathrm{Ker}(p_1+f'^*\,p_2), $
\item $\mathrm{Ker} (\cup [D_5])=\mathrm{Ker}(p_2),$
\end{itemize}

The kernel of 
$$\cup [D_5] : H^2(A\times \mathbb{P}^N,\reel)\subset H^2(X,\reel)\ra H^4(X, \reel)$$
is $$\mathrm{Ker}(q_1+i^*q_2),$$
where the inclusion $H^2(A\times \mathbb{P}^N,\reel)\subset H^2(X,\reel)$ is given by pullback by the composite of $\tau$ with the projection of $A\times A\times \mathbb{P}^N$ on the two last factors.
\end{lem}
Obviously, the lemma remains true after letting any automorphism of $\comp$ act.

\begin{proof}

Let us first prove the assertion about the images by considering the general situation of a blow-up $\tau : \widetilde{Y}\ra Y$ of a complex smooth projective variety along a smooth, but not necessarily irreducible, subvariety $Z$, of codimension everywhere at least $2$. Let $E$ be the exceptional divisor. It is a projective bundle over $Z$. Let $j_Z$ and $j_E$ be the inclusions of $Z$ and $E$ in $Y$ and $\widetilde{Y}$ respectively. It is known (see \cite{V1}, 7.3.3) that there exists a homomorphism $\phi : H^*(Z, \reel)\ra H^*(E, \reel)$, given by excision and the Thom isomorphism, such that the cohomology of $\widetilde{Y}$ is the quotient in the following exact sequence of (non-graded) vector spaces
$$\xymatrix{0\ar[r] & H^*(Z,\reel)\ar[r] & H^*(Y, \reel)\oplus H^*(E, \reel)\ar[r] & H^*(\widetilde{Y}, \reel)\ar[r] & 0}, $$
where the first map is $(j_{Z*}, \phi)$ and the second one is $\tau^*\oplus j_{E*}$

Now let $Z_1, \ldots, Z_n$ be the irreducible components of $Z$, $E_1, \ldots, E_n$ the corresponding irreducible components of $E$, $\tau_{E_i}$ the restriction of $\tau$ to $E_i$, and $j_{Z_i}$ and $j_{E_i}$ the obvious inclusions. Let $x$ be a degree $2$ cohomology class in $Y$. We want to show that the images of the homomorphisms $\cup [E_i]\circ\tau^*$ and $\tau^*\circ\,\cup x$, restricted to degree $1$ cohomology classes in $Y$, are in direct sum. Indeed, since $\tau^*$ is injective on cohomology because $\tau$ is birational, and since the images of $\cup h$ and $\cup \alpha$ in $H^3(A\times A\times\mathbb{P}^N, \reel)$ are in direct sum, as the K\"unneth formula shows, this will prove the assertion.

We have 
$$\cup [E_i]\circ\tau^*=j_{E_i*}\circ\tau_{E_i}^*\circ j_{Z_i}^*,$$
so it is enough to prove that the images of the $j_{E_i*}:H^1(E_i, \reel)\ra H^3(\widetilde{Y}, \reel)$ and of $\tau^* : H^3(Y, \reel)\ra H^3(\widetilde{Y}, \reel)$ are in direct sum. But the map 
$$H^3(Y, \reel)\oplus H^1(E, \reel)\ra H^3(\widetilde{Y}, \reel)$$
has zero kernel for degree reasons, as the exact sequence above shows. This proves the assertion about the images.

Let us now consider the kernels. To compute the first two, it is enough to work on $A\times A\times \mathbb{P}^N$ since $\tau^*$ is injective on cohomology. The K\"unneth formula shows that $h$ has nonzero cup-product with any nonzero element of $H^*(X, \reel)$ coming from $A\times A$.

Consider now cup-product with $\alpha$. Since the real cohomology algebra of an abelian variety is the exterior algebra on the first real cohomology space, the assertion concerning Ker$(\cup \alpha)$ boils down to the following lemma.

\begin{lem}
Let $V$ be a finite dimensional vector space, and let $\alpha$ be a nonzero element of $\bigwedge^2 V$. The kernel of the homomorphism 
$$\wedge\alpha : V\ra \bigwedge^3 V$$
is at most $2$-dimensional.
\end{lem}

\begin{proof}
Let us choose a basis $e_1, \ldots, e_n$ for $V$. The space $\bigwedge^2 V$ has a basis consisting of all the $e_i\wedge e_j$ with $i<j$. Without loss of generality, we can assume that $\alpha$ has a nonzero component on $e_1\wedge e_2$ with respect to this basis. It is then clear that the elements $\alpha\wedge e_3, \ldots, \alpha\wedge e_n$ of $\bigwedge^3 V$ are linearly independent. This shows that the homomorphism $\wedge\alpha : V\ra \bigwedge^3 V$ has rank at least $n-2$, and concludes the proof.
\end{proof}

As for the computation of the other kernels, since the cohomology of a smooth variety is embedded in the cohomology of any smooth blow-up of it, the result is a straightforward consequence of the following general computation.

Consider the situation of two smooth complex projective varieties $B$ and $C$, together with a morphism $f:B\ra C$. Let $\tau : \widetilde{B\times C}\ra B\times C$ be the blow-up of $B\times C$ along the graph $\Gamma$ of $f$. Let $E$ be the exceptional divisor, and let $\tau_E$ be the restriction of $\tau$ to $E$. Let $j_{\Gamma}$ and $j_E$ be the inclusions of $\Gamma$ and $E$ into $\widetilde{B\times C}$ and $B\times C$ respectively. The map 
$$\cup [E] \circ \tau^* : H^*(B\times C,\reel)\ra H^{*+2}(\widetilde{B\times C}, \reel)$$
is equal to 
$$j_{E*}\circ\tau_{E}^*\circ j_{\Gamma}^*.$$
It follows from \cite{V1}, 7.3.3 that $j_{E*}\circ\tau_{E}^*$ is injective, which means that the kernel of $\cup [E] \circ \tau^*$ is equal to the kernel of $j_{\Gamma}^*$. Now the morphism from $B$ to $B\times C$ with coordinates $\mathrm{Id}_B$ and $f$ identifies $\Gamma$ with $B$, and $j_{\Gamma} : \Gamma \hookrightarrow B\times C$ with 
$$(\mathrm{Id}_B\times f) : B\ra B\times C.$$
As a consequence, the kernel of $j_{\Gamma}^*$ is equal to the kernel of the homomorphism 
$$H^*(B\times C,\reel)=H^*(B, \reel)\otimes H^*(C, \reel)\ra H^*(B, \reel)$$
which sends an element of the form $\alpha\otimes\beta$ to $\alpha\cup f^*(\beta)$.

If $b$ and $c$ are complex points of $B$ and $C$, let $p$ and $q$ be the projections from $H^*(B\times C,\reel)$ to $H^*(B, \reel)$ and $H^*(C,\reel)$ induced by the immersions $B\hookrightarrow B\times C, x\mapsto (x,c)$ and $C\hookrightarrow B\times C, x\mapsto (b,x)$. By the K\"unneth formula, if $\gamma$ is a degree $1$ cohomology class in $H^*(B\times C,\reel)=H^*(B, \reel)\otimes H^*(C, \reel)$, we have 
$$\gamma=p(\gamma)\otimes 1 + 1\otimes q(\gamma),$$
which shows that $j_{\Gamma}^*(\gamma)=0$ if and only if $p(\gamma)+ f^*\,q(\gamma)=0$.

This proves that the kernel of 
$$\cup [E] \circ \tau^* : H^1(B\times C,\reel)\ra H^{3}(\widetilde{B\times C}, \reel)$$
is 
$$\mathrm{Ker}(p+f^*q).$$
It is straightforward to check that this equality remains true for degree $2$ cohomology classes in case $B$ or $C$ has no degree $1$ cohomology, since we have then $H^2(B\times C, \reel)=H^2(B, \reel)\oplus H^2(C, \reel)$.

\end{proof}

\subsubsection{Proof of proposition \ref{morph}}

\begin{proof}
Without loss of generality, we can suppose that $\sigma'$ is the identity. Let $\sigma$ be an automorphism of $\comp$ and $\gamma$ be an isomorphism from $H^*(X, \reel)$ to $H^*(X^{\sigma}, \reel)$. We will show that $\gamma$ sends $\reel[D_k]$ to $\reel[D_k^{\sigma}]$ for any $k$ between $1$ and $5$, and $\reel h$ to $\reel h^{\sigma}$. This will be achieved step by step.

\paragraph{The Albanese map.}
We use the same argument as in \cite{V}. Recall that $\pi$ is the natural map from $X$ to $A\times A$. Pullback by $\pi$ gives an isomorphism between $H^1(A\times A, \reel)$ and $H^1(X, \reel)$. The injection 
$$\pi^* : H^*(A\times A, \reel)=\bigwedge H^1(X, \reel) \hookrightarrow H^*(X, \reel)$$
is given by this isomorphism and cup-product. As a consequence, cup-product alone allows us to recover $H^*(A\times A, \reel)$ as a subalgebra of $H^*(X, \reel)$ : it is the algebra $\bigwedge H^1(X, \reel)$. This being true also after letting $\sigma$ act, we get the following.

\begin{lem}\label{ab}
The isomorphism $\gamma$ sends the subalgebra $H^*(A\times A, \reel)$ of $H^*(X, \reel)$ to the subalgebra $H^*(A^{\sigma}\times A^{\sigma}, \reel)$ of $H^*(X^{\sigma}, \reel)$.
\end{lem}

\bigskip

\paragraph{Image of the cohomology classes of the exceptional divisors.}

\begin{lem}
There exists a permutation $\phi$ of $\{1,\dots,5\}$ such that for each $k$ between $1$ and $5$, $\gamma$ sends the line $\reel[D_k]$ to the line $\reel[D_{\phi(k)}^{\sigma}]$.
\end{lem} 

\begin{proof}

From the K\"unneth formula and the computation of the cohomology of a blow-up, we get a splitting
$$H^2(X, \rat)=\pi^* H^2(A\times A, \rat)\oplus \rat h \oplus \bigoplus_{k=1}^5 \rat [D_k].$$
Let $k$ be an integer between $1$ and $5$. The isomorphism $\gamma$ sends $[D_k]$ to some element of $H^2(X^{\sigma}, \reel)$ of the form 
$$\gamma([D_k])=\alpha^{\sigma} + \mu_1 [D_1^{\sigma}]+\ldots+\mu_5 [D_5^{\sigma}]+\nu h^{\sigma},$$
with $\alpha^{\sigma}$ coming from $H^2(A^{\sigma}\times A^{\sigma}, \reel)$. 

We now use lemma \ref{cup}. The map 
$$\cup [D_k^{\sigma}] : H^1(X^{\sigma}, \reel) \ra H^3(X^{\sigma}, \reel)$$
has a $2\dim{A}=4$-dimensional kernel.
Furthermore, the kernel of 
$$\cup\alpha^{\sigma} + \mu_1 [D_1^{\sigma}]+\ldots+\mu_5 [D_5^{\sigma}]+\nu h^{\sigma} : H^1(X^{\sigma}, \reel) \ra H^3(X^{\sigma}, \reel)$$
is the intersection of the kernels of the $\cup \mu_i\,[D_i^{\sigma}]$, $\cup\alpha^{\sigma}$ and $\cup\nu h^{\sigma}$, because the images of these homomorphisms are in direct sum. Since $\cup h^{\sigma}$ is injective on degree $1$ cohomology, we get $\nu=0$. Furthermore, the kernel of $\cup\alpha^{\sigma} : H^1(X^{\sigma}, \reel) \ra H^3(X^{\sigma}, \reel)$ is at most $2$-dimensional unless $\alpha^{\sigma}=0$. This proves $\alpha^{\sigma}=0$, and since the kernels of the $\cup \mu_i\,[D_i^{\sigma}]$ are pairwise distinct $4$-dimensional vector spaces, which implies that the intersection of two of them has dimension at most $3$, this implies that $[D_k]$ is sent to some $\mu_i [D_i^{\sigma}]$. 
This proves the lemma.
\end{proof}

\begin{lem}
The permutation $\phi$ is the identity.
\end{lem}

\begin{proof}
For $k$ between $1$ and $5$, let $F_k$ be the subspace of $H=H^1(X, \reel)$ consisting of elements $\alpha$ such that $\alpha\cup[D_k]=0$. For $\sigma$ an automorphism of $\comp$, let $F_k^{\sigma}$ be the corresponding subspace of $H^{\sigma}=H^1(X^{\sigma}, \reel)$. We determined those spaces in lemma \ref{cup}. From there, and from the actual definition of $f$ and $f'$, one sees that $F_1$ is the only one of the $F_k$ that has a nonzero intersection with two other ones, namely $F_3$ and $F_4$. The same is true for the $F_k^{\sigma}$. As a consequence, $\phi(1)=1$ and $\phi(\{3,4\})=\{3,4\}$, thus $\phi(\{2,5\})=\{2,5\}$. 

\medskip

For $k$ and $k'$ between $1$ and $5$ such that $F_k\cap F_{k'}=0$, let $p_{kk'}$ (resp. $p_{kk'}^{\sigma}$) be the projection along $F_k$ onto $F_{k'}$ (resp. along $F_k^{\sigma}$ onto $F_{k'}^{\sigma}$). Since $\gamma$ sends $\reel[D_k]$ to $\reel[D_{\phi(k)}^{\sigma}]$, it sends $F_k$ to $F_{\phi(k)}^{\sigma}$. The projection $p_{kk'}$ is therefore conjugate to $p_{\phi(k)\phi(k')}^{\sigma}$.

Recall that we chose $f$ and $f'$ to be generators of the complex multiplication of $A$. Let $f^*$ and $f'^*$ be the homomorphisms they induce on the first cohomology group of $A$. Recall that $f^*$ (resp. $f'^*$) has eigenvalues $0$, $a$ and $-a$ (resp. $0$, $a'$ and $-a'$), with $a$ (resp. $a'$) being a generator of $k$ (resp. $k'$). Direct computation shows that, identifying $F_1$ with $H^1(A, \reel)$, we have the equality
$$(p_{21}\circ p_{53})_{|F_1}=1-f^*.$$
As a consequence, the endomorphism $(p_{\phi(2)1}^{\sigma}\circ p_{\phi(5)\phi(3)}^{\sigma})_{|F_1^{\sigma}}$ of $F_1^{\sigma}$ is conjugate to $1-f^{\sigma*}$. 

If $\phi(2)=2$ and $\phi(5)=5$, this imposes $\phi=$Id, as $f^*$ and $f'^*$ have different eigenvalues. 
If  $\phi(2)=5$, $\phi(5)=2$ and $\phi(3)=3$, then 
$$(p_{\phi(2)1}^{\sigma}\circ p_{\phi(5)\phi(3)}^{\sigma})_{|F_1^{\sigma}}=(p_{51}^{\sigma}\circ p_{23}^{\sigma})_{|F_1^{\sigma}}=(1-f^{\sigma *})^{-1},$$
where we identified $F_1^{\sigma}$ with $H^1(A^{\sigma}, \reel)$. Again, consideration of the eigenvalues proves that $(1-f^{\sigma *})^{-1}$ is not conjugate to $1-f^{\sigma*}$. Indeed, $\frac{1}{1-a}$ can't be equal to either $1-a$ or $1+a$, so this case cannot happen. Similarly, we cannot have $\phi(2)=5$, $\phi(5)=2$ and $\phi(3)=4$. This proves that $\phi$ is the identity.
\end{proof}
 
\bigskip

\paragraph{Image of h.}

The only thing left to show is that $\gamma$ sends the line $\reel h$ to the line $\reel h^{\sigma}$. 

\bigskip

Since $\gamma$ is an isomorphism, and because of the preceding paragraph, it sends $h$ to some nonzero multiple of an element of $H^2(X^{\sigma}, \reel)$ of the form 
$$h^{\sigma}+ \alpha^{\sigma} + \lambda_1 [D_1^{\sigma}]+\ldots+\lambda_5 [D_5^{\sigma}],$$
where $\alpha^{\sigma}$ is the pull-back of a class in $H^2(A^{\sigma}\times A^{\sigma}, \reel)$.

The $Z_k$ are pairwise disjoint, so the cup-product of any two different $[D_k]$ is $0$. Furthermore, if $H$ is a generic hyperplane of $\mathbb{P}^N$, then $A\times A\times H$ is disjoint from $Z_1, Z_2, Z_3$ and $Z_4$, which proves that $h\cup[D_k]=0$ for $k$ between $1$ and $4$. This is of course true after conjugation by $\sigma$.

Let $k$ be between $1$ and $4$. Since $[D_k]$ is sent to a nonzero multiple of $[D_k^{\sigma}]$, we conclude that 
$$[D_k^{\sigma}]\cup\alpha^{\sigma} + \lambda_k[D_k^{\sigma}]^2=0$$
for $k$ between $1$ and $4$. We will compute this more explicitely to show that it implies that $\lambda_k$ and $\alpha^{\sigma}$ are both zero.

\bigskip

Let us write $\alpha^{\sigma}=\tau^{\sigma*}\beta^{\sigma}$, where $\tau$ is the blowing-down map and $\beta^{\sigma}$ is a class in $H^2(A^{\sigma}\times A^{\sigma}\times\mathbb{P}^N, \reel)$ coming from $H^2(A^{\sigma}\times A^{\sigma}, \reel)$.
For $k$ between $1$ and $5$, let $h_k\in H^2(D_k, \rat)$ be the first Chern class of the normal bundle of $D_k$ in $X$. It follows from \cite{V1}, lemma 7.32 that the cohomology of $D_k$ is a free module on the cohomology of $Z_k$, with basis $1, h_k, \ldots, h_k^{N+1}$, since the codimension of $Z_k$ in $A\times A \times\mathbb{P}^N$ is $N+2$. For simplicity, we will drop the $\sigma$ superscript when applied to morphisms.

The self-intersection formula gives, for any positive integer $a$, 
$$[D_k^{\sigma}]^a=\tilde{j}_{k*}((h_k^{\sigma})^{a-1}),$$
where $\tilde{j}_{k}$ is the inclusion of $D_k$ in $X$. As a consequence, we can compute, for any positive $a$ and nonnegative $b$, 
$$[D_k^{\sigma}]^a\cup (\tau^*\beta^{\sigma})^b=\tilde{j}_{k*}((h_k^{\sigma})^{a-1}\cup \tau_k^* j_k^* (\beta^{\sigma})^b).$$
In particular, for $k$ between $1$ and $4$, we get
$$\tilde{j}_{k*}(\tau_k^* j_k^* (\beta^{\sigma}) + \lambda_k h_k^{\sigma} )=0.$$

Using \cite{V1}, 7.3.3 as before, we see  that the map 
$$\widetilde{j}_{k*}: H^{2}(D_k^{\sigma}, \reel)\ra H^{4}(X^{\sigma}, \reel)$$
is injective. This proves that 
$$\tau_k^* j_k^* (\beta^{\sigma}) + \lambda_k h_k^{\sigma}=0$$
in $H^2(D_k^{\sigma}, \reel)$, for any $k$ between $1$ and $4$, which is equivalent to $\lambda_k=0$ and $j_k^* (\beta^{\sigma})=0$. Making $k=1$ and $k=2$, then using lemma \ref{cup}, this proves that $\beta^{\sigma}$ is in the kernel of $p_1$ and $p_2$, hence is zero.

\bigskip

We thus have shown that $\gamma$ sends $h$ to some nonzero multiple of an element of $H^2(X^{\sigma}, \reel)$ of the form 
$$h^{\sigma}+ \lambda [D_5^{\sigma}].$$
We want to show that $\lambda$ is zero. The next lemma concludes the proof.

\begin{lem}
Let $\lambda$ be a real number. Then
$$(h+\lambda [D_5])^{N+1}=0 \Leftrightarrow \lambda=0.$$
The same is true after conjugation by $\sigma$.
\end{lem}

\begin{proof}
We obviously have $h^{N+1}=0$, so let us suppose $(h+\lambda [D_5])^{N+1}=0.$
Let $H$ be a hyperplane in $\mathbb{P}^N$. Recall that $h$ is the pullback by the blowing-down map $\tau$ of the cohomology class $g$ of $A\times A\times H$ in $A\times A\times\mathbb{P}^N$. Using the preceding computation, we get 
$$\widetilde{j}_{5*}(\sum_{i=0}^{N}\binom{N+1}{i+1}\lambda^{i+1} h_5^{i}\cup \tau_5^*j_5^*(g^{N-i}))=0.$$
Using \cite{V1}, 7.3.3 again, and noticing that the map $$j_{5*} : H^*(Z_5, \reel)\ra H^*(A\times A\times\mathbb{P}^N, \reel)$$
is injective, we see that the homomorphism 
$$\widetilde{j}_{5*}: H^{2N}(D_5, \rat)\ra H^{2N+2}(X, \rat)$$
is injective. We thus get
$$\sum_{i=0}^{N}\binom{N+1}{i+1}\lambda^{i+1} h_5^{i}\cup \tau_5^*j_5^*(g^{N-i})=0,$$
hence, since the $h_5^{i}$ are linearly independent over $H^*(Z_5, \reel)$ for $i< \mathrm{codim}\, Z_5=N+2$, we obtain 
$$\lambda^{i+1} \tau_5^*j_5^*(g^{N-i})=0$$
for $i$ between $0$ and $N$. For $i=N$, this proves that $\lambda=0$.
\end{proof}

Now, since $h^{N+1}=0$, we have 
$$(h^{\sigma}+ \lambda [D_5^{\sigma}])^{N+1}=0.$$
The preceding lemma thus shows that $\lambda=0$, which concludes the proof.
\end{proof}

\subsubsection{Proof of proposition \ref{recover}}

\begin{proof}

Once we know lemma \ref{cup}, the proof of this proposition is purely formal. Indeed, let $\sigma$ be an automorphism of $\comp$. We want to recover, given the abstract graded algebra $H^*(X^{\sigma}, \reel)$ together with the lines $\reel[D_1^{\sigma}], \ldots, \reel[D_5^{\sigma}]
$ and $\reel h^{\sigma},$ the space $H^1(A^{\sigma}, \reel)$ with the action of $f^{\sigma}$ and $f'^{\sigma}$, the space $H^2(\mathbb{P}^N, \reel)$ and the restriction map 
$$i^{\sigma*} : H^2(\mathbb{P}^N, \reel)\ra H^2(A^{\sigma}, \reel)$$
induced by the projective embedding of $A^{\sigma}$.

We are going to use the restriction maps $p_1^{\sigma}$ and $p_2^{\sigma}$ from $H^1(X^{\sigma}, \reel)=H^1(A^{\sigma}\times A^{\sigma}, \reel)$ to $H^1(A^{\sigma}, \reel)$ induced by the inclusions of the first and the second factor respectively. From lemma \ref{cup}, we see that the data we have allow us to recover, in the space $H^1(X^{\sigma}, \reel)$, the subspaces 
\begin{displaymath}
\mathrm{Ker}(p_1^{\sigma}),\, \mathrm{Ker}(p_2^{\sigma}), \,\mathrm{Ker}(p_1^{\sigma}+p_2^{\sigma}), \,\mathrm{Ker}(p_1^{\sigma}+f^{\sigma*}p_2^{\sigma}), \,\mathrm{Ker}(p_1^{\sigma}+f'^{\sigma*}p_2^{\sigma}).
\end{displaymath}

Giving those subspaces is equivalent to giving the vector space $H^1(A^{\sigma}, \reel)$ together with the action of $f^{\sigma}$ and $f'^{\sigma}$ on it. Indeed, the first two subspaces determine a splitting of $H^1(X^{\sigma},\reel)$ into two subspaces isomorphic to $H^1(A^{\sigma}, \reel)$. The third one is the graph of a specific isomorphism between them, which allow us to identify them -- actually, using the opposite of this particular isomorphism in order to get the right sign. The last two subspaces are then the graphs of the endomorphisms $-f^{\sigma}$ and $-f'^{\sigma}$ of $H^1(A^{\sigma}, \reel)$. We thus recover the real vector space $H^1(A^{\sigma}, \reel)$ with the action of $f^{\sigma}$ and $f'^{\sigma}$.

The same procedure allows us to recover the other data through the second part of lemma \ref{cup}. Indeed, the line $\reel h$ is equal to the space $H^2(\mathbb{P}^N, \reel)\subset H^2(X^{\sigma}, \reel)$. In the previous paragraph, we showed how to recover the subspace $\mathrm{Ker}(p_1^{\sigma})=H^1(A^{\sigma}, \reel)\subset H^1(X^{\sigma}, \reel)$, which allows us to recover $$H^2(A^{\sigma}, \reel)=\bigwedge^2 H^1(A^{\sigma}, \reel)\subset H^2(X^{\sigma}, \reel),$$
this being the image of the cohomology of $A^{\sigma}$ under the pull-back by the projection on the second factor. We thus obtained the subspace $H^2(A^{\sigma}\times \mathbb{P}^N, \reel)\subset H^2(X^{\sigma}, \reel)$ and its direct sum decomposition 
$$H^2(A^{\sigma}\times \mathbb{P}^N, \reel)=H^2(A^{\sigma}, \reel)\oplus H^2(\mathbb{P}^N, \reel).$$
Using lemma \ref{cup}, we obtain the graph of the opposite of the homomorphism
$$i^{\sigma *} : H^2(\mathbb{P}^N, \reel) \ra H^2(A^{\sigma}, \reel).$$
\end{proof}

{\par\vskip.3truein\relax
\leftline{\hskip2truein\relax Fran\c{c}ois Charles}
\leftline{\hskip2truein\relax \'Ecole Normale Sup\'erieure}
\leftline{\hskip2truein\relax 45 rue d'Ulm}
\leftline{\hskip2truein\relax 75\thinspace 230 PARIS CEDEX 05}
\leftline{\hskip2truein\relax FRANCE}
\leftline{\hskip2truein\relax\tt francois.charles@ens.fr}
\relax}

\end{document}